\documentclass[proceedings]{stacs}
\stacsheading{2009}{565--576}{Freiburg}
\firstpageno{565}

\usepackage[T1]{fontenc}
\usepackage{amsmath}
\usepackage{amsthm}
\usepackage{amssymb}

\newcommand{\rank}              {\mathit{rank}}

\newcommand{\dom}{\mathit{dom\, }}
\newcommand{\tra }{\mathit{Tr}}
\newcommand{\comp}{\mathit{Comp}}

\begin{document}
\title{On the Borel inseparability of game tree languages}
\author[lab1]{S. Hummel}{Szczepan Hummel}
\author[lab1]{H. Michalewski}{Henryk Michalewski}
\author[lab1]{D. Niwi\'nski}{Damian Niwi\'nski}
\thanks{All authors are supported by Grant~N206 008 32/0810.}
\address[lab1]{Faculty of Mathematics, Informatics, and Mechanics\newline Warsaw University, Poland}
\email{{shummel,henrykm,niwinski}@mimuw.edu.pl} 


\begin{abstract}
The game tree languages can be viewed as an automata-theoretic counterpart of parity games on graphs.
They witness the strictness of the index hierarchy of alternating tree automata, as well as the 
fixed-point hierarchy over binary trees.

We consider a game tree language of the first non-trivial level, where Eve 
can force that 0 repeats from some moment on, and its dual, where Adam
can  force that 1  repeats from some moment on. Both these sets (which 
amount to one up to an obvious renaming) are complete in the class of co-analytic sets. We show that they cannot be
separated by any Borel set, hence {\em a fortiori\/} by any
weakly definable set of trees.

This settles a case left open  by L.Santocanale and A.Arnold,
who have thoroughly investigated the separation property within the $\mu $-calculus
and the automata index hierarchies.
They showed that separability fails in general for non-deterministic automata 
of type $\Sigma^{\mu }_{n} $, starting from level 
$n=3$, while our result settles the missing case $n=2$.
\end{abstract}

\maketitle

\section*{Introduction}
In 1970 Rabin~\cite{Rabin70} proved 
the following 
property:
If a set of infinite trees can be defined both by an  {\em existential\/}
and by a  {\em universal\/} 
sentence of monadic second order logic  then it can 
also be defined in a weaker logic, with quantification restricted to {\em finite\/} sets.
An automata-theoretic counterpart of this fact~\cite{Rabin70,Saoudi}
states 
that if a tree language,  as well as  its complement, are both recognizable by B\"uchi automata
(called {\em special\/} in~\cite{Rabin70}) then they are also recognizable by weak
alternating automata. Yet another formulation, in terms of the
$\mu $-calculus~\cite{an_weak}, states that if a  tree language 
is definable both by a $\Pi^{\mu }_2 $-term (i.e.,  with a  pattern
$\nu \mu $) and a $\Sigma^{\mu }_2 $-term ($\mu  \nu $),
then it 
is also definable by an alternation free term, i.e., one in $\comp (\Pi^{\mu }_1 \cup  \Sigma^{\mu }_1)$.
This last formulation 
gives rise to a question 
if the equation
\[
\Pi^{\mu }_{n} \cap \Sigma^{\mu }_{n} =
\comp (\Pi^{\mu }_{n-1} \cup  \Sigma^{\mu }_{n-1})
\]
holds  on all levels of the fixed-point  hierarchy.
Santocanale and Arnold showed~\cite{arnold_santo}, rather surprisingly, that it is
not the case for $n \geq 3 $. They exhibit  a series of 
``ambiguous''  properties, expressible by terms in $\Pi^{\mu }_{n} $ and 
in $\Sigma^{\mu }_{n} $, but not in $\comp (\Pi^{\mu }_{n-1} \cup  \Sigma^{\mu }_{n-1}) $.
On positive side however, they discover  a more subtle generalization of Rabin's result,
which continues to hold on the higher stages of  the  hierarchy.
\smallskip

Let us explain it  at a more abstract level, with $\mathcal L $ (``large'') and $\mathcal S $ (``small'')
being  two classes of subsets of some  universe $U$. Consider the following properties.
\begin{quote}
{\em Simplification.} Whenever  $L$ and its  complement $\bar{L} $ are both in $\mathcal L $,
they are also in $\mathcal S $. 
\end{quote}
\begin{quote}
{\em Separation.} Any two disjoint sets $L,M \in {\mathcal L} $ are separated by some set $K$ in $\mathcal S $
(i.e., $L \subseteq K \subseteq U-M $).
\end{quote}
Note that (given some $\mathcal L $ and $\mathcal S $) separation implies simplification, but in general not
{\em vice versa\/}. In topology, 
it is well known (see, e.g., \cite{Kechris94}) that the separation property holds for
\begin{quote}
${\mathcal L} = $  analytic ($ {\bf \Sigma}^1_1 $) subsets of a  Polish space (e.g., $\{ 0,1 \}^{\omega } $),\\
${\mathcal S} = $  Borel sets,
\end{quote}
 but fails for ${\mathcal L} = $ co-analytic sets ($ {\bf \Pi}^1_1$)
and ${\mathcal S}$ as above.
On the other  hand both classes enjoy the 
simplification property (which amounts to the Suslin Theorem).
\smallskip

In this setting, Rabin's result establishes the simplification property for 
\begin{quote}
${\mathcal L} = $ 
B\"uchi
definable tree languages ($\Pi^{\mu }_{2} $ in the fixed-point hierarchy),\\
${\mathcal S} = $  weakly definable tree languages
($\comp (\Pi^{\mu }_{1} \cup  \Sigma^{\mu }_{1})$).
\end{quote}
 A closer
look at the original proof  reveals that a (stronger)
separation property also holds for these classes.
\smallskip

Santocanale and Arnold~\cite{arnold_santo} showed in  turn that 
 the separation property holds 
for 
\begin{quote}
${\mathcal L} = $ tree languages  recognizable by 
{\em non-deterministic\/}  automata of level $\Pi^{\mu }_{n} $,\\
${\mathcal S} = $ 
tree languages definable by  fixed-point terms in
$\comp (\Pi^{\mu }_{n-1} \cup  \Sigma^{\mu }_{n-1})$,
\end{quote}
for the remaining case of $n \geq 3$.
On the negative side,
they showed that  the separation property fails 
for ${\mathcal L} $ consisting of tree languages recognizable by 
{\em non-deterministic\/}  automata of level $\Sigma^{\mu }_{n} $, for
$n \geq 3 $, leaving open the case of $n=2$. In fact, their proof reveals that, 
in the case under consideration,
even a (weaker) simplification property fails
(see~\cite{arnold_santo}, section 2.2.3).  
As for  $\Sigma^{\mu }_{2} $ however, 
the simplification property does hold,
because of  Rabin's result\footnote{\label{stopka}If a set and its complement
are recognized by non-deterministic co-B\"uchi automata then they
are also both recognized by alternating B\"uchi automata~\cite{an_Buchi}, and hence 
by non-deterministic B\"uchi automata, and hence are weakly
definable~\cite{Rabin70}.}.  
For this reason, the argument of Santocanale and Arnold cannot be extended to
the class $\Sigma^{\mu }_{2} $.
In the present paper, we show that the separation property
fails also in this case,
completing the missing point in the classification of~\cite{arnold_santo}.
\smallskip

We use a topological argument and show in fact a somewhat stronger result, exhibiting 
two disjoint languages recognized by non-deterministic tree automata with 
co-B\"uchi condition (i.e., $\Sigma^{\mu }_{2} $),
which cannot be separated by any {\em Borel\/} set (in a standard Cantor-like topology
on trees). The languages in question are 
the so-called {\em game tree languages\/} (of level (0,1)), which were 
used in~\cite{bradfield:1998} 
(and later also in~\cite{A:hierarchy})
in the proof of the strictness
of the fixed-point hierarchy over binary trees. More specifically, one of
these languages consists of the trees labeled in $\{ 0,1 \} \times \{ \exists , \forall \} $,
such that in the induced game (see definition below) {\em Eve\/} has a strategy to force 
only 0's from some moment on. The second is the twin copy of the first and  consists of
those trees that {\em Adam\/} has a strategy to force 
only 1's from some moment on. 
\smallskip

The wording introduced above differs slightly from the standard terminology of
descriptive set theory, where a {\em separation property\/} of a class $\mathcal L$ means
our  property with
${\mathcal S} = \{ X : X, \bar{X} \in {\mathcal L} \} $ 
(see~\cite{Kechris94}). To emphasize the distinction, following~\cite{Addi},
we will refer to the latter as to the {\em first  separation property\/}.
In this setting, the first separation property holds for the class
of B\"uchi
recognizable tree languages, but it fails for the co-B\"uchi languages, similarly
as it is  the case of the analytic {\em vs.} co-analytic sets, mentioned above.
This may be read as an evidence of  a strong analogy between the B\"uchi class
and $ {\bf \Sigma}^1_1 $. In fact, Rabin~\cite{Rabin70} early observed
that the  B\"uchi tree languages are definable by existential sentences of monadic
logic, and hence
analytic. We show however that, maybe surprisingly, the converse is not true,
by exhibiting  an  analytic  tree language, recognized by
a parity (Rabin) automaton, but not by 
any  B\"uchi automaton.
\smallskip

{\em Note.} The fixed-point hierarchy discussed above provides an obvious  context of
our  results, but in the paper we do not rely on the $\mu $-calculus concepts or methods.
For definitions of relevant concepts, we refer an interested reader to the
work by Santocanale and Arnold~\cite{arnold_santo} or, e.g., 
to~\cite{mubook}.

\section{Basic concepts} \label{prelimy}

Throughout the paper, $\omega $ stands for the set of natural numbers. 

\paragraph{Metrics  on trees}
A full binary tree over a finite alphabet $\Sigma  $ (or shortly a tree,
if confusion does not arise) is represented as a mapping
$t : \{ 1,2 \}^* \to \Sigma  $. 

We consider the classical topology {\em \`a la Cantor\/} on $T_{\Sigma }$
induced by the metric
\begin{equation} \label{distance}
d( t_1, t_2 ) = \left\{\begin{array}{ll}
0
&\mbox{ if $t_1 = t_2 $}\\
2^{-n} \mbox{ with } n = \mbox{min} \{ |w| : t_1 (w) \neq t_2 (w) \}
& \mbox{ otherwise }
\end{array}\right.
\end{equation}
It is well-known and easy to see that if
$\Sigma $ has at least two elements then
$T_{\Sigma }$ with this topology is
homeomorphic to the Cantor discontinuum $\{ 0,1 \}^{\omega }$.
Indeed, it is enough to fix a  bijection $\alpha : \omega \to \{ 1,2 \}^*$,
and a mapping (code) $C: \Sigma \to \{ 0, 1 \}^* $, such that
$C(\Sigma )$ forms a maximal antichain w.r.t.~the prefix ordering.
Then $T_{\Sigma } \ni t \mapsto C \circ  t \circ \alpha \in \{ 0,1 \}^{\omega }$
is a desired homeomorphism.
We assume that the reader is familiar with the basic concepts of 
set-theoretic topology
(see, e.g., \cite{Kechris94}).
The {\em Borel sets\/} over $T_{\Sigma }$ constitute the least family containing 
open sets and closed under  complement and countable union. 
The  {\em Borel relations\/} are defined similarly, starting with open relations
(i.e., open subsets of $T_{\Sigma }^n$, for some $n$, considered with 
product topology). The {\em analytic\/} (or $ {\bf \Sigma}^1_1 $) sets
are those representable by
\[
L = \{ t :  (\exists t')\, R (t, t')   \}
\]
where $R \subseteq T_{\Sigma } \times T_{\Sigma } $ 
is a Borel relation. The  {\em co-analytic\/}  (or $ {\bf \Pi}^1_1 $)
sets are the complements of analytic sets. A continuous mapping 
$f : T_{\Sigma } \to  T_{\Sigma } $ {\em reduces\/} a tree language
$A \subseteq  T_{\Sigma }$ to $B  \subseteq  T_{\Sigma } $ if
$f^{-1} (B) = A $. As in complexity theory, a set $L \in {\mathcal K}$ is 
{\em complete\/} in class $\mathcal K $ if 
all sets in this class reduce to it.

\paragraph{Non-deterministic automata.}
A non-deterministic tree automaton over trees in $T_{\Sigma } $
with a 
parity acceptance condition\footnote{Currently 
most frequently used in the literature, these automata are well-known to be
equivalent to historically previous 
automata with the Muller or Rabin conditions~\cite{thomas:hb-fl}.}
is presented as
${A}\, =\, \langle\Sigma, Q,q_{I},
\tra , \rank \rangle $, where  $Q$ is
a finite set of {\em states} with an {\em initial state} $q_{I}$,
$\tra  \subseteq Q\times \Sigma \times Q  \times Q $ is a set of
{\em transitions},  and $ \rank : Q \rightarrow \omega $ is the 
{\em ranking} function.    
 A transition $(q,\sigma  ,p_1, p_2)$ is usually written
$q \stackrel{\sigma }{\to } p_1, p_2 $.

A {\em run\/} of $ A$ on a tree $t \in T_{\Sigma}$ is itself a 
$Q$--valued tree $\rho  : \{ 1,2\}^* \rightarrow Q$  such that 
$\rho (\varepsilon  ) = q_I$, and, for each $w \in \dom (\rho )$,
$\rho (w) \stackrel{t(w)}{\to } \rho (w1), \rho (w2) $ is
a transition in $\tra $.
A {\em path} $P = p_0 p_1 \ldots \in \{ 1,2  \}^{\omega }$
in $\rho $ is {\em accepting}
if the highest rank occurring infinitely often along it
is even, i.e.,
$\limsup_{n \rightarrow \infty }{\rank} ( \rho (p_0 p_1 \ldots p_n ) )$  is { even\/}. 
A \emph{run is accepting} if so are all its paths.
A tree language $T(A)$ {\em recognized} by $A$ consists of those trees
in $T_{\Sigma }$ which admit an accepting run.

The {\em Rabin--Mostowski index\/} of an automaton $A$ is the pair 
$(\mbox{min} (\rank (Q)),\linebreak \mbox{max} (\rank (Q)))$; without loss of
generality, we may assume that $\mbox{min} (\rank (Q)) \in \{ 0, 1\} $.

An automaton with the {\em  Rabin--Mostowski index\/} $(1,2) $ is called
a  {\em B\"uchi automaton\/}. Note that a B\"uchi automaton accepts a tree $t$
if, on each path, some state of rank 2 occurs infinitely often.  We refer to
the tree languages recognizable by B\"uchi automata as to 
 {\em B\"uchi (tree) languages\/}. The  {\em co-B\"uchi languages\/} are the
complements of  B\"uchi languages.  It is known that if a tree language is
recognized by a non-deterministic automaton of index $(0,1) $ then it is
co-B\"uchi\footnote{It follows, in particular, from the equivalence of the
non-deterministic and alternating B\"uchi automata~\cite{an_Buchi},
mentioned in footnote~\ref{stopka}.}; the converse is not true in
general (see the languages $M_{i,k} $ in Example~\ref{przyklad_Rabina} below).

\begin{example} \label{przyklad_Rabina}
Let
\[
L = \{ t \in T_{\{ 0,1 \}} :
(\exists P) \, \limsup_{n \rightarrow \infty } t (p_0 p_1 \ldots p_n) = 1 \}
\]
This set is recognized by a B\"uchi automaton with transitions
\[
\begin{array}{lllllll}
q/p  \stackrel{0}{\to } q,T; \;  & \; &
q/p \stackrel{1}{\to } p,T; \;  & \; &
T \stackrel{(0/1)}{\to } T,T; \; & \; &
\\
q/p  \stackrel{0}{\to } T,q; \;  & \; &
q/p \stackrel{1}{\to } T,p; \;  & \; &
& &
\end{array}
\] 
with $\rank (q) = 1$ and  $\rank (p) = \rank (T) = 2$.
Rabin~\cite{Rabin70} showed that its complement $\bar{L} $ cannot be recognized by any
B\"uchi automaton, but it is recognizable by an (even deterministic) automaton of index $(0,1) $
\[
\begin{array}{lllllll}
0/1  \stackrel{0}{\to } 0,0 ; \;  & \; & 
0/1  \stackrel{1}{\to } 1,1 ; \;  & \; & \rank (i) = i ; \;  & \; \mathit{for}  & i = 0,1. 
\end{array}
\] 
This last set can be generalized to the so-called parity languages  (with $i \in \{ 0,1\} $)
\[
M_{i,k} = \{ t \in T_{\{ i, \ldots , k \}} :
(\forall P) \, 
\limsup_{n \rightarrow \infty } t (p_0 p_1 \ldots p_n) \; \;  \mathit{is} \; \;  \mathit{  even}
\}
\]
which are all co-B\"uchi but require arbitrary high indices~\cite{Niw86}.
It can also be showed that all languages $M_{i,k} $ (except for 
$(i,k) = (0,0),(1,1),(1,2) $) are complete in the class of co-analytic sets 
$ {\bf \Pi}^1_1 $ (see, e.g., \cite{gapa}).
\end{example}

The class of languages which are simultaneously B\"uchi and co-B\"uchi has numerous
characterizations mentioned in the introduction; all 
these characterizations  easily imply that such
sets are Borel (even of finite Borel rank).  

\begin{example} \label{przyklad_Pawla}
Consider the set $\bar{L} = M_{0,1} $ of Example~\ref{przyklad_Rabina}, and its twin copy obtained
by the renaming $0 \leftrightarrow 1 $,
\[
M_{0,1}' = \{ t \in T_{\{ 0,1 \}} : (\forall P) \, 
\liminf_{n \rightarrow \infty } t (p_0 p_1 \ldots p_n) = 1 \}.
\]
The sets $M_{0,1}$ and $M_{0,1}'$ are disjoint, co-B\"uchi and, as we have already noted, 
 $ {\bf \Pi }^1_1 $ complete.
They can be separated by a set $K$ of trees\footnote{This argument is due to
Pawe{\l } Milewski.}, such that on the {\em rightmost\/}
branch, there are only finitely many 1's
\[  
K = \{ t \in T_{\{ 0,1 \}} : \limsup_{n \rightarrow \infty } t (\underbrace{22\ldots 2}_n) = 0 \}
\]
(i.e., $M_{0,1} \subseteq K \subseteq T_{\{ 0,1 \}} - M_{0,1}' $).
The set $K$ can be presented as a countable union of closed sets
\[  
K = \bigcup_m  \; \{ t \in T_{\{ 0,1 \}} :  (\forall n \geq m)\, 
t (\underbrace{22\ldots 2}_n) = 0 \}
\]
so it is on the level ${\bf \Sigma}^0_2 $ (i.e., $F_{\sigma } $) of the Borel hierarchy.
The membership in the Borel hierarchy can also be seen trough an automata-theoretic argument by showing that $K$ is simultaneously B\"uchi and co-B\"uchi. Indeed it can be recognized by
an (even deterministic) automaton with co-B\"uchi condition
\[
\begin{array}{lllllllll}
0/1  \stackrel{0}{\to } T,0 ; \;  & \; & 
0/1  \stackrel{1}{\to } T,1 ; \;  & \; & 
T \stackrel{(0/1)}{\to } T,T; & \; &
\rank (i) = i , \;   \mathit{for} \;  \;   i = 0,1, & \rank (T) = 0 , 
\end{array}
\] 
as well as by a (non-deterministic) B\"uchi automaton
\[
\begin{array}{llllllll}
q  \stackrel{(0/1)}{\to } T,q / p ; \;  & \; & 
p  \stackrel{0}{\to } T,p ; \;  & \; & 
T \stackrel{(0/1)}{\to } T,T; & \; &
\rank (q) = 1, \;   \; & 
\rank (p) = \rank (T) = 2.
\end{array}
\] 
\end{example}
We will see in the next section that a Borel separation of
 co-B\"uchi languages is not always possible.

\section{Inseparable pair} \label{para}
Let
\[
\Sigma = \{ \exists , \forall \} \times  \{ 0, 1\} ,
\]
we denote by $\pi_i $ the projection on the {\em i\/}th component
of $\Sigma  $.
With each $t \in T_{\Sigma }$, we associate a game $G(t) $, played by two players, 
{\em Eve\/} and {\em Adam\/}. The positions of Eve are those nodes $v$,
for which  $\pi_1 (t(v)) = \exists $, the remaining nodes are positions of
Adam. For each position $v$, it is possible to move to one of its successors,
$ v1$ or $v2$.
The players  start in the root and then 
move  down the tree, 
thus forming an infinite path 
$P= (p_0 p_1 p_2 \ldots )   $. 
The successor is selected by Eve or Adam depending on who is the
owner of the position
$p_0 p_1 \ldots p_{n-1}  $.   The play is  won  by Eve if
\[
\limsup_{n \to   \infty } \pi_2 \left(
t (p_0 p_1 \ldots p_{n}) \right) = 0
\]
i.e., 1 occurs only finitely often, otherwise Adam is the winner.
A strategy for Eve  selects a move for each of
her positions; it is winning if any play consistent with the strategy
is won by Eve. We say that Eve {\em wins\/} the game $G(t) $ if she has a winning strategy.
The analogous concepts  for Adam are defined  similarly.

A reader familiar with the {\em parity games\/} (\cite{emerson:jutla:91}, see
also \cite{thomas:hb-fl}) has noticed of course that 
the games $G(t)$  are a special case of these
(with the index $(0,1)$).
\smallskip

Now let
\[
W_{0,1} = \{ t : \mathit{Eve} \; \; \mathit{wins} \; \; G(t)  \}
\]
We also define a set 
$W_{0,1}' \subseteq T_{\Sigma } - W_{0,1} $, consisting of those trees $t$,
where Adam has a strategy which guarantees 
him not only to win in $G(t) $, but
also to force a stronger condition, namely 
\[
\liminf_{n \to   \infty } \pi_2 \left(
t (p_0 p_1 \ldots p_{n}) \right) = 1 .
\]
It should be clear that $W_{0,1}' $ can be obtained from  $W_{0,1} $ by applying
(independently on each component) a renaming  $0 \leftrightarrow 1 $,
 $\exists \leftrightarrow \forall $. Thus, the sets  $W_{0,1} $ and  $W_{0,1}' $ are
disjoint,  but have identical topological and automata-theoretic properties.
\smallskip

Let us see that the set   $W_{0,1} $ can be recognized by a non-deterministic automaton
of index $(0,1) $; it is enough to take the states  $\{ 0,1  \} \cup \{ T \} $,
 with $\rank (T) = 0 $, and $\rank (\ell ) = \ell $, for $\ell \in  \{ 0,1  \}$,  the initial 
state $0$,
and transitions
\[
\begin{array}{lllllll}
\ell \stackrel{( \forall , m)}{\to } m,m; \;  & \; &
\ell \stackrel{(\exists , m  )}{\to } m,T; \;  & \; &
\ell \stackrel{( \exists ,m )}{\to } T,m; \; & \; &
T \stackrel{( Q , m  )}{\to } T,T,
\end{array}
\] 
with $m \in \{ 0,1 \} $, and $Q \in \{ \exists , \forall \} $.
Hence, the sets  $W_{0,1} $ and  $W_{0,1}' $ are co-B\"uchi
(c.f. the remark before Example~\ref{przyklad_Rabina}).

We are ready to state the main result of this paper.
\begin{theorem} \label{glowne_tw}
The  sets $W_{0,1} $ and  $W_{0,1}' $ cannot be separated by any Borel set.
\end{theorem}
\proof
The proof relies on the following.
\begin{lemma} \label{glowny_l}
For any Borel set $B \subseteq T_{\Sigma }$,  there is a continuous function
$f_B :  T_{\Sigma }\to  T_{\Sigma } $, such that 
\[
\begin{array}{lll}
u \in B  & \Rightarrow & f_B (u) \in W_{0,1} \\
u \not\in B & \Rightarrow & f_B (u) \in W_{0,1}'
\end{array}
\]
\end{lemma}

\proof
Note that $f_B$ is required to reduce simultaneously $B$ to  $W_{0,1} $ and $T_{\Sigma } - B$
 to  $W_{0,1}' $. We proceed by induction on the complexity of the set $B$.

Note first that if $B$ is clopen (simultaneously closed and open) then it is enough to fix two
trees $t \in W_{0,1} $ and  $t' \in W_{0,1}' $, and define  $f_B$ by 
\[
\begin{array}{lll}
u \in B  & \Rightarrow & f_B (u) = t \\
u \not\in B & \Rightarrow & f_B (u) = t'
\end{array}
\]
Also note that, by symmetry of  the sets $W_{0,1} $ and   $W_{0,1}' $, the claim
for $B$ readily implies the claim for the complement $T_{\Sigma } - B $. (Specifically,
$f_{B'} $ is obtained by composing $f_{B} $ with a suitable renaming.)

Finally note that the space $T_{\Sigma } \approx \{ 0,1 \}^{\omega }$
has a countable basis consisting of clopen sets.

\begin{figure}[tp]
  \centering
\includegraphics[width=0.3\textwidth]{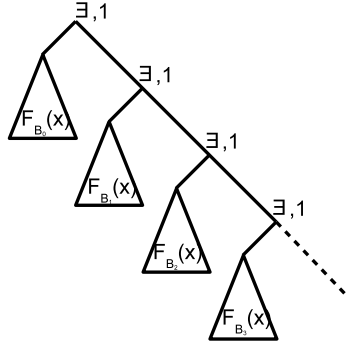}
  \caption{Induction step for $\bigcup_n B_n $. \label{figura}}
 \end{figure}

Then, in order to complete the proof,
it remains to settle the induction step for $B = \bigcup_{n < \omega } B_n $.
Assume that we have already the reductions $f_{B_n} $ satisfying the claim, for $n < \omega $.
Given $u \in T_{\Sigma }  $, we construct a tree $f_B (u) $, by labeling
the rightmost path by $( \exists , 1) $, and letting a subtree in the node $2^n 1 $
be  $f_{B_n} (u)$
(see Figure~\ref{figura}).
In symbols,
\[
\begin{array}{lllll}
 f_B (u) (2^n) & = & ( \exists , 1) &  &  \\
 f_B (u) (2^n 1 v) & = &  f_{B_n} (u) (v),  & \mathit{for} \; \; n < \omega ,
&  v \in \{ 1,2 \}^* \, .
\end{array}
\]
Since all the functions $f_{B_n} $ are continuous, the resulting $f_B $ is
continuous as well.
Now, if $u \in B_{m} $, for some $m$, then Eve has an obvious winning strategy:
follow the rightmost path and turn left in $2^m $, then use the winning strategy
on the subtree $f_{B_m} (u) $, which exists, by induction hypothesis.
 
If, however, $(\forall n)\, u \not\in B_n $ then Adam can win the game with the
stronger winning criterion, required in the definition of $W_{0,1}' $.
Indeed, he can do so as soon as Eve enters any of the subtrees  $f_{B_n} (u) $
(by induction hypothesis), but he also wins if Eve remains forever on the
rightmost path.

This proves the claim for  $f_B$, and thus completes the proof of the lemma.
\qed
We are ready to complete the proof of the theorem. Suppose that there is a Borel
set $C$, such that $W_{0,1} \subseteq C  \subseteq  T_{\Sigma } - W_{0,1}' $. The
claim of the lemma immediately implies that
\[
\begin{array}{lll}
u \in B  & \Rightarrow & f_B (u) \in C \\
u \not\in B & \Rightarrow & f_B (u) \in T_{\Sigma } - C.
\end{array}
\]
Thus any Borel set $B$ over $T_{\Sigma }$ is reducible to $C$, but this
is clearly impossible, as it would contradict the strictness of
the Borel rank hierarchy in the Cantor discontinuum $\{ 0,1 \}^{\omega } $ 
(see, e.g., \cite{Kechris94}).
\qed

Since the sets $W_{0,1} $ and  $W_{0,1}' $ are recognizable by non-deterministic automata
of index $(0,1) $, Theorem~\ref{glowne_tw} settles the case of $n=2 $, missing in Section
2.2.3 of~\cite{arnold_santo}, devoted to the failure of separation property for
non-deterministic automata of type $\Sigma^{\mu }_n $ and the class
$\comp (\Pi^{\mu }_{n-1} \cup  \Sigma^{\mu }_{n-1})$.

In the terminology introduced at
the end of introduction, we can state the following.

\begin{corollary}
The class of co-B\"uchi tree languages does not have the first separation property.\qed
\end{corollary}

This may be contrasted  with the positive result of~\cite{Rabin70}. As we have mentioned
in the introduction, Rabin's original proof essentially shows this property
for the class of B\"uchi tree languages,
although it is not explicitly stated there. For the sake of completeness, we
sketch the argument below, following closely the 
$\mu $-calculus version of~\cite{an_weak} (based on
the original proof of~\cite{Rabin70}).

\begin{theorem}[Rabin]
The class of B\"uchi tree languages has the first separation property.
\end{theorem}
\proof
Let $A$ and $B$ be two non-deterministic B\"uchi automata,
such that $T(A) \cap T(B) = \emptyset $. 
We will refer to the states of rank 2  as
to {\em accepting\/} states (of the corresponding automaton).
A {\em cut\/} (of a tree)  is a finite maximal
antichain in $\{ 1,2 \}^* $ with respect to the prefix
ordering $\leq $.
For two cuts $X,Y$ we let $Y > X $ if $Y$ lies
below  $X$, i.e.,
$(\forall y \in Y) \, (\exists x \in X) \, y > x $. 
It is easy to see that a run $\rho $ of  a  B\"uchi automaton is accepting if,
for each cut $X$, there is a cut  $Y > X $, labeled by the accepting states
(i.e., $(\forall y \in Y) \, \rank (\rho (y)) = 2 $). 
We inductively define a sequence of tree languages $K^n_q $, for each 
state $q$ of $A$, and $n \geq 0$.

The set  $K^0_q $ consists of all trees $t$ 
which admit some run (not necessarily accepting) of $A$  
  starting from $q$
($q$-run, for short).  The set  $K^{n+1}_q $ comprises those trees  $t$,  which
admit a $q$-run $\rho $, such  that, for each cut $X$, there exists
a cut $X' > X $, and a run $\rho'$,  with the following properties:
\begin{itemize}
\item   $\rho'$  agrees with $\rho $ until the cut $X$,
\item all states in $\rho' (X') $ are accepting,
\item $(\forall v \in X')\, $, the subtree of $t$ 
rooted in $v$  (in symbols $t.v $)  belongs to $K^{n}_p $, where $p = \rho' (v) $.
\end{itemize}
It follows by induction on $n$ that $T(A) \subseteq K^n_{q_I}   $, where
$q_I $ is the initial state of $A$.
Now let $n_A $ and $n_B$ be the numbers of states of
$A$ and $B$, respectively, and let $M = 2^{n_A \cdot n_B} + 1 $. We claim
that $K^{M}_{q_I} $ separates $T(A)$ and  $T(B)$. 
We already know that $T(A) \subseteq K^{M}_{q_I} $. For the sake of contradiction, 
suppose that 
$t \in K^{M}_{q_I} \cap T(B) $, and let $\rho' $ be an
accepting run of $B$ on $t$. 

Using the inductive definition of $K^{M}_{q_I} $, we can 
construct a sequence of cuts
$X_1 < X'_1 < \ldots < X_M < X'_M $, and a run $\rho $ of $A$ on $t$,
such that
\begin{itemize}
\item  $(\forall i \leq M)$ all states in $\rho (X_i) $ are accepting,
\item $(\forall i \leq M, \, \forall v \in X_i)\, t.v \in K^{M-i}_{\rho (v)} $, 
\item  $(\forall i \leq M)$ all states in $\rho' (X'_i) $ are accepting.
\end{itemize}
By the choice of $M$, there exist $1 \leq k < \ell \leq M $, such that
\[
\{ (\rho (u), \rho' (u) ) : u \in X_k \}   = 
\{ (\rho (v), \rho' (v) ) : v \in X_{\ell } \}
\]
Note that, by construction,
\[
 X_k < X'_k < X_l
\]
with all states in $\rho' (X'_k ) $ accepting.
Hence, by a standard tree-pumping argument, we can construct a new tree along with
two accepting runs: by $A$ and by $B$, contradicting 
$T(A) \cap T(B) = \emptyset $.

It remains to show that the language $K^{M}_{q_I}$ is both
B\"uchi and co-B\"uchi. A direct construction of two B\"uchi automata
would be somewhat cumbersome, but one can use here any of the
characterizations of this intersection class  mentioned above.
In the proof given in~\cite{an_weak}, it is shown that the sets  $K^{n}_{q}$
are definable in the alternation-free $\mu $-calculus. A reader familiar
with monadic second-order logic can easily see that
these languages are definable in its {\em weak\/} fragment, i.e., with
quantifiers restricted to finite sets. This is enough as well, according
to the characterization given by Rabin~\cite{Rabin70}.
\qed

\section{Broken analogy}
A reader familiar with descriptive set theory may think of another
inseparable pair of recognizable tree languages, induced by a classical
example (\cite{Kechris94}, section 33.A).
We will explain why it would not be useful for our purpose.
Let us now  consider  non-labeled trees, i.e., subsets  $T \subseteq \omega^* $
closed under initial segments.  
They
can be viewed  as elements of the Cantor discontinuum $\{ 0,1 \}^{\omega } $ by 
fixing a
bijection $\iota : \omega \to \omega^* $ and identifying a tree $T$ with
its characteristic function, given by $f_T (n) = 1 $ iff $\iota (n) \in T $.
In particular, we can discuss  topological properties of sets
of such trees.
As before, $P \in \omega^{\omega } $ is
a path in a tree $T$ if all finite prefixes of $P$ are in $T$.
Let
\begin{eqnarray*}
\mbox{WF} & = & \{ T : \mbox{$T$ has no infinite path } \} \\
\mbox{UB}  & = & \{ T : 
\mbox{$T$ has {\em exactly one\/} infinite path } \}
\end{eqnarray*}
Both sets are known to be $ {\bf \Pi}^1_1 $-complete, although
the membership of UB in  $ {\bf \Pi}^1_1 $ is not obvious, and is
the subject of one of Lusin's theorems (Theorem 18.11 in~\cite{Kechris94}).
WF and UB  are also known to be inseparable by Borel 
sets  (\cite{Kechris94}, section  35, see also \cite{Becker}). 
Now, it is not difficult to ``encode'' these sets as languages of
labeled binary trees, which turn out to be recognizable by parity
automata.
In~\cite{gapa} a continuous reduction of WF to $M_{0,1} $ was used to
show that the latter set is complete in $ {\bf \Pi}^1_1 $
(Example~\ref{przyklad_Rabina} above).  Let
\[
\begin{array}{llll}
\mbox{UB}_{\mathit{bin}} & = & 
\{ t \in T_{\{ 0, 1 \} } : & \mbox{ there is {\bf exactly one} path $P$} \\
   & &  & \mbox{with } 
\limsup_{n \rightarrow \infty } t (p_0 p_1 \ldots p_n ) = 1  \}
\end{array}
\]
It is easy to construct a non-deterministic automaton accepting 
this language;
one can also  assure that this
automaton is non-ambiguous, i.e., for each accepted tree,  has exactly
one accepting run. From considerations above, one can deduce 
that the sets $T_{0,1} $ and $\mbox{UB}_{\mathit{bin}} $ are
inseparable by Borel languages. However,
the language  $\mbox{UB}_{\mathit{bin}} $ is not co-B\"uchi.

\begin{proposition}\label{kontr_przykladzik}
The language $\overline{\mbox{\em UB}_{\mathit{bin}}} $ is
recognizable and analytic, but not B\"uchi.
\end{proposition}
\proof
Let us call a path with infinitely many 1's {\em bad\/}. So the
above language consists of trees that have either none or at least two
bad paths.
Rabin~\cite{Rabin70}  shows that the language
$T_{0,1} $ (no bad paths) cannot be recognized by a B\"uchi
automaton, by constructing  a correct tree which by
pumping argument can be transformed to a tree with 
{\em exactly
one\/} bad path (mistakingly accepted by the hypothetical automaton).
So this classical argument applies to the language 
$\mbox{UB}_{\mathit{bin}} $ without any changes.
\qed
As we have argued in the introduction, this example somehow breaks
the analogy between the class of B\"uchi recognizable tree languages and
that of analytic sets.  It turns out that the topological 
complexity, and the
automata-theoretic complexity, although closely related, 
do not always coincide.

\section{Conclusion}
The automata-theoretic hierarchies, in particular the index hierarchies for
non-de\-ter\-mi\-nis\-tic and alternating tree automata, are studied because of
the issues of expressibility and complexity. Typically, the higher  the level
in the hierarchy, the  higher  the expressive power of automata, but also 
the complexity
of the related algorithmic problems (like emptiness or inclusion).
Once the strictness of the hierarchy is established~\cite{bradfield:1997,bradfield:1998},
the next important problem  is an {\em effective simplification\/},
i.e., determining the exact level of an object (e.g., a tree language)
in the hierarchy. The problem is generally unsolved (see~\cite{thomas_chris:2008}
for a recent development in this direction). One may expect that a better understanding of
structural properties of the hierarchy can bring a progress also in this problem.
We believe that ideas coming from  descriptive set theory,
like separation and reduction properties, uniformization, or
completeness, can be helpful here.
\smallskip

The inseparable pairs of co-analytic sets are common in mathematics.
Natural examples include 
the 
set of all continuous real--valued functions on the unit interval $[0,1]$ which 
are everywhere differentiable together with the set of all continuous real--valued 
functions on the unit interval $[0,1]$ which are 
not differentiable in exactly one point, 
but as in this case, other examples usually reflect the same pattern 
of WF {\em vs.\/} UB (c.f. \cite{Becker}). 
In contrast, our pair presented in Section~\ref{para} is very symmetric:
the two sets are  copies of each other up to a symbolic renaming.
Recently,  Saint Raymond~\cite{Raymond} 
established  that the pair WF {\em vs.\/} UB 
is complete (in the sense of Wadge)
with respect to all coanalytic pairs in the Cantor set.
In
the proof 
he uses an interesting example of another  complete coanalytic pair, 
which exhibits certain symmetric properties. Building on his results, in subsequent
work, we show that the  pair $W_{0,1} $, $W_{0,1}' $, has an analogous
completeness property.
\bigskip

Our example shows that the first separation property fails for the 
co-B\"uchi class ($\Sigma^{\mu }_2 $ in the fixed-point hierarchy)
while, by Rabin results~\cite{Rabin70}, it holds for the 
B\"uchi class ($\Pi^{\mu }_2 $).
By this we have also settled a missing case in a classification by
Santocanale and Arnold~\cite{arnold_santo}. However, these authors
were  interested in the {\em relative\/} separation property
(as explained in our introduction),  as they primarily wanted to find
if the {\em ambiguous class\/} $ \Pi^{\mu }_{n} \cap \Sigma^{\mu }_{n}$ can
be effectively captured by $\comp (\Pi^{\mu }_{n-1} \cup  \Sigma^{\mu }_{n-1}) $.
As this coincidence turned out to fail for $n \geq 3 $, it is meaningful to ask if
the status of the first separation property
established  for the B\"uchi/co-B\"uchi classes, continues to hold
for the
higher-level classes $\Pi^{\mu }_n $/$\Sigma^{\mu }_n $. That is,
if two disjoint sets definable  in $\Pi^{\mu }_n $ can always be separated
by a set in $ \Pi^{\mu }_{n} \cap \Sigma^{\mu }_{n}$.
(A similar question for $\Sigma^{\mu }_n $, with expected answer negative.) 
In our opinion, it is an interesting 
problem, which may challenge for a better understanding of the
topological structure of recognizable languages
above  $ {\bf \Pi}^1_1 \cup {\bf \Sigma }^1_1 $.

\newpage\null
\end{document}